\documentclass[12pt,a4paper]{article}
\usepackage{amssymb}
\usepackage{amsmath}
\usepackage{amsthm}
\usepackage{xypic}

\theoremstyle{plain}
\newtheorem{theorem}{Theorem}
\newtheorem{proposition}[theorem]{Proposition}

\theoremstyle{remark}
\newtheorem{remark}[theorem]{Remark}
\theoremstyle{definition}
\newtheorem{definition}[theorem]{Definition}
\newtheorem{example}[theorem]{Example}
\newtheorem{lemma}[theorem]{Lemma}
\newtheorem{conjecture}[theorem]{Conjecture}

\def\varinjlim_#1{\lim\limits_{\longrightarrow\atop{#1}}}

\def\End{\mathop{\rm End}\nolimits}

\def\Aut{\mathop{\rm Aut}\nolimits}
\def\Hom{\mathop{\rm Hom}\nolimits}
\def\id{\mathop{\rm id}\nolimits}
\def\BSU{\mathop{\rm BSU}\nolimits}
\def\BU{\mathop{\rm BU}\nolimits}
\def\SU{\mathop{\rm SU}\nolimits}

\def\PGL{\mathop{\rm PGL}\nolimits}

\def\Z{\mathop{\rm Z}\nolimits}
\def\KSU{\mathop{\rm KSU}\nolimits}
\def\AB{\mathop{\rm AB}\nolimits}
\def\K{\mathop{\rm K}\nolimits}
\def\diag{\mathop{\rm diag}\nolimits}
\def\pt{\mathop{\rm pt}\nolimits}

\begin{document}

\author{A. V. Ershov}
\title{Floating bundles and their applications}
\date{}
\maketitle
\section {The homotopic functor, connected with floating
algebra bundles}

The aim of this section is to define the homotopic
functor to category of Abelian groups, connected
with the special classes of bundles
with fiber matrix algebra
or projective space. The results, which will be represented
(without proofs\footnote{see some
proofs in the articles (in Russian)
[1] Ershov A. V. O gomotopicheskikh svoistvakh rassloenii
so strukturnoi gruppoi avtomorfizmov matrichnikh algebr // Vestnik
Moskovskogo universiteta, ser. Matematika. Mekhanika. ¹6, 1999,
s. 56-58, [2] Ershov A. V.
O $\K$-teorii rassloenii na matrichnie algebri // UMN,
T. 55, Vip. 2, 2000, c. 137-138.})
below, were obtained in the author's dissertation.

Let $X$ be a finite $CW$-complex, $\widetilde M_n:=X\times
M_n(\mathbb{C})$ be the product bundle over $X$ with fiber
$M_n(\mathbb{C})$ (where $M_n(\mathbb{C})$ is the algebra of all
$n\times n$ matrices over $\mathbb{C}$).

The next definition is motivated by well-known result that every
vector bundle over compact base is the subbundle in a product
bundle.

\begin{definition}
\label{bd}
Let $A_k$ be a locally trivial
bundle over $X$ with fiber $M_k(\mathbb{C})$.
Assume that there exists a bundle map
 \begin{equation}
 \nonumber
 \begin{array}{c}
 \diagram
 A_k\rrto^\mu \drto && \widetilde M_{kl}\dlto \\
 &X&
 \enddiagram
 \end{array}
 \end{equation}
such that for any point $x\in X$
the fiber $(A_k)_x\cong M_k(\mathbb{C})$
is embedded (by the restriction $\mu_x$ of $\mu$)
into the fiber $({\widetilde M}_{kl})_x\cong
M_{kl}(\mathbb{C})$ as a central simple subalgebra.
Then the triple $(A_k,\mu ,\widetilde M_{kl})$ is called the {\it algebra
bundle (AB)} over $X$. Moreover if $k$ and $l$ are
relatively prime (i. e. if their greatest common divisor $(k,l)=1$), then
the triple $(A_k,\mu ,\widetilde M_{kl})$ is called the {\it floating
algebra bundle (FAB)}.

\end{definition}

There exists the similar class of bundles with fiber projective
space.

\begin{definition}
\label{bdd}
Let $X$ be a finite $CW$-complex; $P^{k-1}$ and
$Q^{l-1}$ be two locally trivial bundles (over $X$) with fibers
$\mathbb{C} P^{k-1}$ and $\mathbb{C} P^{l-1}$ respectively.
By $\widetilde {\mathbb{C} P}^{kl-1}$
denote the product bundle over $X$
with fiber $\mathbb{C} P^{kl-1}.$
Assume that there exists a bundle map
\begin{equation}
\nonumber
\begin{array}{c}
\diagram
\stackrel{ }{\stackrel{ }{P^{k-1}{\mathop {\times} \limits_X}
Q^{l-1}}}\rrto ^\lambda  \drto &&
\stackrel{{\displaystyle \widetilde {{\mathbb{C}} P}^{kl-1}}}{ } \dlto \\
&X&
\enddiagram
\end{array}
\end{equation}
(where $P^{k-1}{\mathop {\times}\limits_X} Q^{l-1}$ is
the product of the bundles over $X$)
such that for any point $x\in X$ the fiber
$(P^{k-1}{\mathop{\times}\limits_X}
Q^{l-1})_x \cong \mathbb{C} P^{k-1}\times \mathbb{C} P^{l-1}$
is embedded into the fiber $(\widetilde {\mathbb{C} P}^{kl-1})_x\cong
\mathbb{C} P^{kl-1}$ by the map of Segre. Then the bundle
$P^{k-1}{\mathop{\times}\limits_X}Q^{l-1}$
is called the {\it bundle of Segre's product
(BSP)}, and under the condition $(k,l)=1$ the
{\it floating bundle of Segre's product (FBSP)}.
\end{definition}

Note that the projective space $\mathbb{C}P^{n-1}$ and
the matrix algebra $M_n(\mathbb{C})$ have the same group
of automorphisms $\PGL(n).$
There exists the special class of principal bundles
such that any bundle satisfies the conditions of
Definition \ref{bd} or \ref{bdd}
is associated with a bundle from this class.
Therefore we get one-to-one correspondence between AB and BSP
(respectively FAB and FBSP) over $X$.

A morphism from AB $(A_k,\mu ,\widetilde M_{kl})$
to AB $(B_m,\nu ,\widetilde M_{mn})$ (over $X$)
is a pair $(f,g),$ where $f\colon A_k\hookrightarrow B_m$
and $g\colon \widetilde M_{kl}\hookrightarrow \widetilde M_{mn}$
are bundle maps such that their restricts to any fiber are
monomorphisms of central simple algebras and the corresponding
quadratic diagram
\begin{equation}
\begin{array}{ccc}
\widetilde{M}_{kl} & \stackrel{g}{\hookrightarrow} & \widetilde{M}_{mn}\\
\scriptstyle{\mu} \uparrow && \uparrow \scriptstyle{\nu}\\
A_k & \stackrel{f}{\hookrightarrow} & B_m
\end{array}
\end{equation}
is commutative. Note that the morphism $(f,g)$ exists only if
$k\mid m, l\mid n.$

Morphisms of BSP may be defined analogously.

The set of all subalgebras in $M_{kl}(\mathbb{C})$,
that are isomorphic to $M_k(\mathbb{C}),$ is parametrized by
the homogeneous manifold $Gr'_{k,kl}.$ This manifold
may be considered as an analog of Grassmannian
manifold.
Moreover, the set of Segre's maps $\mathbb{C}P^{k-1}
\times \mathbb{C}P^{l-1}\hookrightarrow \mathbb{C}P^{kl-1}$
(where $\mathbb{C}P^{kl-1}$ is fixed)
is parametrized by the same manifold $Gr'_{k,kl}.$

Note that the stable (i. e. $\pi_r(Gr'_{k,kl})$
under the condition $r<2\min\{k,l\}$)
homotopic groups of $Gr'_{k,kl}$ are
the same as for $\BSU.$

Let us consider the canonical AB
$({\cal A}'_k,\mu' ,\widetilde{\cal M}'_{kl})$
over $Gr'_{k,kl},$ defined in
the following way. The fiber of the bundle ${\cal A}'_k$ over
a point $x\in Gr'_{k,kl}$ is the subalgebra in $M_{kl}(\mathbb{C}),$
which corresponds to $x$. Thus the bundle ${\cal A}'_k$
and its embedding into the product bundle
$\widetilde{\cal M}'_{kl}$ are uniquely defined.
If $(k,l)=1$, then $({\cal A}'_k,\mu'
,\widetilde{\cal M}'_{kl})$ is FAB.

The canonical BSP (if $(k,l)=1$, then FBSP)
${\cal P}^{'k-1}{\mathop {\times}\limits_{Gr'_{k,kl}}}{\cal Q}^{'l-1}$
over $Gr'_{k,kl}$ is defined in a similar way.

By the standard method of reduction of a structural group
we can replace the noncompact manifold $Gr'_{k,kl}$
by the homotopy-equivalent compact manifold $Gr_{k,kl}.$
Fibers of the canonical AB
$({\cal A}_k,\mu ,\widetilde{\cal M}_{kl})$
(respectively canonical BSP ${\cal P}^{k-1}{\mathop {\times}
\limits_{Gr_{k,kl}}}{\cal Q}^{l-1}$~\label{spaceBSP})
over $Gr_{k,kl}$ have a concordant Hermitian (respectively
K\"{a}hlerian) structure.

Since the categories of AB and BSP (respectively FAB and FBSP)
are equivalent, we shall formulate the next results only
in terms of AB and FAB.

The following Proposition shows that if $(k,l)=1,$ then
$Gr_{k,kl}$ and the canonical FAB $({\cal A}_k,\mu,{\cal M}_{kl})$
are the classifying space and the universal FAB
(for FAB of the form $(A_k,\mu,\widetilde{M}_{kl})$
over a finite $CW$-complexes) respectively.

\begin{proposition}
Let $X$ be a finite $CW$-complex,
$\dim X<2\min\{k,l\}$,
$(A_k,\mu ,\widetilde M_{kl})$ be a FAB over $X$,
i. e. $(k,l)=1$.
Then there exists a map (unique up to homotopy)
$\varphi \colon X\to Gr_{k,kl}$ such that
$\varphi^* ({\cal A}_k,\mu ,\widetilde {\cal M}_{kl})\cong
(A_k,\mu ,\widetilde M_{kl})$
(i. e. $\varphi$ is a classifying map for
$(A_k,\mu ,\widetilde M_{kl})$).
\end{proposition}

Now let us consider an embedding
$i_{kl,mn}: Gr_{k,kl}\hookrightarrow Gr_{m,mn}$
induced by an embedding of algebras $M_{kl}(\mathbb{C})
\hookrightarrow M_{mn}(\mathbb{C})\;
(\Rightarrow k|m,l|n).$
Let the map $i_{kl,mn*}: \pi_r(Gr_{k,kl})\rightarrow \pi_r(Gr_{m,mn})$
be the homomorphism induced by $i_{kl,mn}.$
Let $\gamma_{k,l}$ be the generator of $\pi_{2s}(Gr_{k,kl}),$ then
$i_{kl,mn*}(\gamma_{k,l})=\frac{(m,n)}{(k,l)}\gamma_{m,n}$
(here we assume that $s<\min\{ k,l\}$; recall
that in this case $\pi_2(Gr_{k,kl})=0,$ $\pi_{2s}(Gr_{k,kl})=
\pi_{2s}(\BSU)=\mathbb{Z},s>1,$
$\pi_{2s+1}(Gr_{k,kl})=\pi_{2s+1}(\BSU)=0$).

It follows from the last result that if $(m,n)=1\;
(\Rightarrow (k,l)=1)$ then the map $i_{kl,mn}\colon
Gr_{k,kl}\rightarrow Gr_{m,mn}$ induces an
isomorphism of the stable homotopic groups.
Therefore the homotopic groups are stabilized
under the passage to the direct limit
$\varinjlim_{(k,l)=1}Gr_{k,kl}.$
This explains the reason of separate study of FAB.
In general case the localization of the homotopic groups
occurs. Note that the theory of FAB (which will be stated below)
cannot be reduced to the theory of usual principal bundles:
embedding $\mu$ plays an important role. Moreover
a bundle $A_k,$ which forms a FAB $(A_k,\mu,\widetilde{M}_{kl}),$
has a special form.

We shall study only FAB below.

It may be proved that for any sequence
$\{ k_j,l_j\}_{j\in \mathbb{N}}$
of pairs, satisfies the following conditions
 $$
 {\rm(i)} \quad k_j,l_j\to\infty;\quad{\rm(ii)}\ k_j|k_{j+1},\
 l_j|l_{j+1};\quad{\rm(iii)}\ (k_j,l_j)=1,
 $$
the spaces $\varinjlim_j Gr_{k_j, k_jl_j}$
are homotopy-equivalent.
By $\varinjlim_{(k,l)=1}Gr_{k,kl}$ or $Gr$
denote this unique homotopic type.

Let us consider a FAB of the form
$(\widetilde M_k,\tau ,\widetilde M_{kl}).$
Let $\tau$ be the map $\tau \colon X\times
M_k(\mathbb{C})\rightarrow X\times M_{kl}(\mathbb{C})$
such that $\tau(x,T)=(x,T\otimes E_l)$ (where $x\in X$,
$T\in M_k(\mathbb{C}),$ $E_l$
is the matrix unit of order $l\times l$, and $T\otimes E_l$
denotes the Kronecker product of matrix). Then the FAB
$(\widetilde{M}_k,\tau,\widetilde{M}_{kl})$ is called {\it trivial}.

\begin{definition}
{\it Isomorphism from FAB} $(A_k,\mu ,\widetilde M_{kl})$ {\it to
FAB}
$(C_k,\nu ,\widetilde M_{kl})$ over $X$ is a pair of bundle maps
$f:A_k\rightarrow C_k$, $g:\widetilde M_{kl}\rightarrow \widetilde
M_{kl}$ such that the following conditions hold:
\begin{itemize}
\item[1)]
for any point $x\in X$ the fiber $(A_k)_x$
(respectively the fiber $(\widetilde{M}_{kl})_x$)
is embeded into the fiber $(C_k)_x$
(respectively into the fiber $(\widetilde{M}_{kl})_x$)
by the restriction of $f$ (respectively of $g$)
as a central simple
subalgebra (in particular $f_x,$ $g_x$ are homomorphisms),
\item[2)]
the following diagram
\begin{equation}
\nonumber
\begin{array}{c}
\diagram \widetilde M_{kl} \rto ^g & \widetilde M_{kl}
\\ A_k\uto ^\mu \rto ^f& C_k\uto _\nu
\enddiagram
\end{array}
\end{equation}
is commutative.
\end{itemize}
\end{definition}

\begin{definition}
\label{floating}
FAB $(A_k,\mu ,\widetilde M_{kl})$ and $(B_m,\nu ,\widetilde
M_{mn})$ are called {\it stable equivalent} if there
exist a sequence of pairs of natural numbers $\{t_i,u_i\}$,
$1\leq i\leq s$
and a corresponding sequence of FAB $(A_{t_i},\mu _i,
\widetilde M_{t_iu_i})$
such that the following conditions hold:
\begin{itemize}
\item[(i)]
$\{t_1,u_1\}=\{k,l\},\{t_s,u_s\}=\{m,n\};$
\item[(ii)]
$(t_it_{i+1},u_iu_{i+1})=1$ for $s>1,\: 1\leq i\leq s-1;$
\item[(1)]
$(A_{t_1},\mu _1,\widetilde M_{t_1u_1})=(A_k,\mu ,\widetilde
M_{kl}),$ $(A_{t_s},\mu _s,\widetilde M_{t_su_s})=(B_m,\nu
,\widetilde M_{mn});$
\item[(2)]
$(A_{t_i},\mu _i,\widetilde M_{t_iu_i})\otimes (\widetilde
M_{t_{i+1}},\tau ,\widetilde M_{t_{i+1}u_{i+1}})\cong (A_{t_{i+1}}
,\mu_{i+1},\widetilde M_{t_{i+1}u_{i+1}})\otimes
(\widetilde M_{t_i},\tau,\widetilde M_{t_iu_i})$
\end{itemize}
for
$s>1,\: 1\leq i\leq s-1,$ where $(\widetilde M_{t_i},\tau ,\widetilde
M_{t_iu_i})$ is the trivial FAB.
\end{definition}

The space $\varinjlim_{(k,l)=1}Gr_{k,kl}$ is the classifying space
for FAB with the just defined relation of the stable equivalence.
This motivates Definition \ref{floating}.

By $[(A_k,\mu,\widetilde{M}_{kl})]$ denote the stable equivalence
class of the FAB $(A_k,\mu,\widetilde{M}_{kl}).$

Now let us show that the set of all
stable equivalence classes of
FAB over $X$
is an Abelian group with respect to the operation,
induced by the tensor product of FAB.
For this we need the following Proposition.

\begin{proposition}
For any pair $\{ k,l\}$ such that $\rm{(i)}$ $(k,l)~=1,$
$\rm{(ii)}$ $2\min\{k,l\}>\dim X$
every equivalence class of FAB over $X$ has a representative
of the form $(A_k,\mu ,\widetilde M_{kl})$.
\end{proposition}

Suppose
$(A_k,\mu ,\widetilde M_{kl})$
and $(B_m,\nu ,\widetilde M_{mn})$ are FAB over~$X$.
If $(km,ln)=1$, then
it is clear that $(A_k\otimes B_m,\mu \otimes \nu ,\widetilde M_{klmn})$
is the FAB. Then by definition, put $[(A_k,\mu,\widetilde{M}_{kl})]\cdot
[(B_m,\nu,\widetilde{M}_{mn})]=[(A_k\otimes B_m,\mu \otimes \nu,
\widetilde{M}_{klmn})].$
Otherwise, applying the previous Proposition,
we can replace $(B_m,\nu ,\widetilde M_{mn})$
by an another representative $(B_{m'},\nu',\widetilde M_{m'n'})$
of the same equivalence class such that $(km',ln')=1$.

It is clear that the product of stable equivalence classes
is well defined. The unit of this operation is the class of
trivial FAB
$[(\widetilde{M}_k,\tau ,\widetilde{M}_{kl})]$.
For given FAB $(A_k,\mu ,\widetilde{M}_{kl})$
consider the subbundle $B_l$ in $\widetilde{M}_{kl}$
such that $(B_l)_x=\Z_{(\widetilde{M}_{kl})_x}((A_k)_x)$ (where
$\Z_B(A)$ is the centralizer of a subalgebra $A$
in an algebra B) for any point $x\in X$.
This also defines the embedding
$\nu \colon B_l\hookrightarrow {\widetilde M}_{kl}.$
Thus the FAB $(B_l,\nu ,\widetilde{M}_{kl})$ is defined.
The inverse element for
$[(A_k,\mu ,\widetilde{M}_{kl})]$ is the class of FAB
$(B_l,\nu,\widetilde{M}_{kl}).$

The next theorem sums the basic results obtained up to now.

\begin{theorem}
The set of all stable equivalence classes of FAB
with respect to the operation, induced by
the tensor product of FAB, defines the contravariant
homotopic functor to the category of Abelian groups.
This functor is denoted by~$\widetilde{\AB}^1$.
Its represented space is the space
$\varinjlim_{(k,l)=1}Gr_{k,kl}$
with the structure of $H$-group, defined by the maps
$Gr_{k,kl}\times Gr_{m,mn}\rightarrow Gr_{km,kmln}$ for $(km,ln)=1$
(the last maps are induced by the tensor product of algebras).
\end{theorem}

\begin{definition}
Let $(A_k,\mu,\widetilde{M}_{kl})$ be an AB.
The bundle $A_k$ (that is considered as
a locally trivial bundle with the structural
group $\Aut M_k(\mathbb{C})\cong \PGL(k)$) is
called the {\it base} of the AB $(A_k,\mu ,\widetilde M_{kl}).$
\end{definition}

In the next Lemma the important property of bundles
(with fiber matrix algebra), which are bases of FAB,
is established.

\begin{lemma}
\label{biglem}
Let $X$ be a finite $CW$-complex, $\dim X<2\min\{k,m\}$.
Then the following conditions are equivalent:
\begin{itemize}
\item[(i)]
$A_k$ is a base of a FAB over $X$;
\item[(ii)]
for any~$m,\; 2m>\dim X,$
there exists a bundle $B_m$ with fiber $M_m(\mathbb{C})$ such that
$A_k\otimes{\widetilde M}_m\cong B_m\otimes {\widetilde M}_k$;
\item[(iii)]
for some $m,$ $(k,m)=1,$ an isomorphism
$A_k\otimes {\widetilde M}_m\cong
B_m\otimes {\widetilde M}_k$ holds.
\end{itemize}
Moreover for any pair of bundles
$(A_k,B_m)$ over $X$ such that
$(k,m)=1$ and $A_k\otimes \widetilde M_m\cong B_m\otimes \widetilde
M_k,$ there exists a unique stable equivalence class of
FAB over $X$ which has representatives of forms
$(A_k,\mu ,\widetilde M_{kn})$ and $(B_m,\nu,\widetilde M_{mn})$
(for any sufficiently large $n$ such that $(km,n)=1$).
This class is denoted by $[(A_k,B_m)].$
\end{lemma}

The most interesting implication in this Lemma is
${\rm (iii)}\Rightarrow {\rm (i)}.$ In its proof,
which is based on the obstruction theory,
the relatively primality of $k$ and $m$ plays
the fundamental role.

Further in this section we describe a connection between
the functor $\widetilde{\AB}^1$ and the reduced $\K$-functor
$\widetilde{\KSU}.$

Let $X$ be a finite $CW$-complex, $\xi _k$ be
a complex vector $\SU(k)$-bundle over $X$ of a rank $k,$
$2k>\dim X.$
By $[n]$ denote the trivial bundle of rank $n$
over $X.$ Take a positive integer m such that $(k,m)=1$
and $m>k.$
Let us consider the pair
$$\xi _k\otimes [k]-[k(k-1)],\; (\xi _k\oplus [m-k])\otimes
[m]-[m(m-1)]$$ of virtual bundles of virtual
dimensions $k$ and $m$ respectively.
The condition $\dim X<2k$ implies that
there exists a unique (up to isomorphism) $k$-dimensional
geometric representative of the stable equivalence
class of $\xi_k\otimes [k]-[k(k-1)].$
By $\xi^\dag_k$ denote this geometric representative.
Let $(\xi _k\oplus [m-k])^\dag$ be the analogous $m$-dimensional
representative of the class of
$(\xi_k\oplus [m-k])\otimes [m]-[m(m-1)].$

The next proposition is almost obviously.

\begin{proposition}
$(\End \xi^\dag_k)\otimes \widetilde M_m\cong (\End (\xi _k
\oplus [m-k])^\dag)\otimes \widetilde M_k.$
\end{proposition}

Let $[\xi_k]$ be the $\widetilde{\KSU}$-equivalence class
of the bundle $\xi_k.$

It follows from Lemma
\ref{biglem} that if $(k,m)=1,$ then
there exists a unique stable equivalence class of FAB
over $X,$ corresponding to the pair $\xi^\dag_k,
(\xi_k\oplus [m-k])^\dag.$
Clearly, that this stable equivalence class of FAB
is independent of a choice
of representative of the equivalence class $[\xi_k].$
Moreover the following Proposition holds.

\begin{proposition}
Let $k,\: m$ be the positive integers such that
$(k,m)=1,$ $2\min \{k,m\}>\dim X.$
Then the map
$$\widetilde{\KSU}(X)\longrightarrow \widetilde{\AB}^1(X),\quad
[\xi_k]\longmapsto [(A_k,B_m)]$$ (where
$A_k=\End \xi_k^\dag,\; B_m=\End(\xi _k\oplus [m-k])^\dag,$
and $[(A_k,B_m)]$ is the stable equivalence class of FAB,
corresponding to the pair $(A_k,B_l)$ as in Lemma~\ref{biglem})
is well defined bijection
of the sets.
\end{proposition}

In particular, any base of FAB has the form
$\End \xi_k^\dag$
for some $\SU(k)$-bundle $\xi_k.$

It follows from the previous Proposition
that the spaces
$\varinjlim_{(k,l)=1}
Gr_{k,kl}$ and $\BSU$ are homotopy-equivalent.

Denote by $\BSU_\otimes$ the space $\BSU$ with the structure
of $H$-group, induced by the tensor product
of virtual $\SU$-bundles of virtual dimension $1$.

\begin{theorem}
The $H$-group $\varinjlim _{(k,l)=1}Gr_{k,kl}$
with respect to the operation, induced by the tensor product of
FAB, is isomorphic to the $\BSU_\otimes$
as an $H$-group.
\end{theorem}

Hence the group $\widetilde{\AB}^1(X)$ is isomorphic to
the multiplicative group (in sense of the theory of formal
groups) of the ring $\widetilde{\KSU}(X),$ i. e.
to the group (since $X$ is the finite $CW$-complex)
$\widetilde{\KSU}(X)$ with respect to the operation
$\xi \circ \eta =\xi +\eta +\xi \eta \; (\xi,\: \eta \in
\widetilde{\KSU}(X)).$

By means of FAB,
the structure of $H$-group on $\BSU_\otimes$
may be described in geometric terms.
Let us consider, for example, the finding
of an inverse element.
Recall that for any FAB $(A_k,\mu ,\widetilde M_{kl})$
there exists the FAB $(B_l,\nu ,\widetilde M_{kl})$
such that for any point $x\in X$ its base $B_l$
has the centralizer $\Z_{(\widetilde{M}_{kl})_x}
(A_k)_x$ as the fiber over $x$
and $\nu$ is the corresponding embedding.
The equivalence class of the FAB $(B_l,\nu ,\widetilde M_{kl})$
is the inverse element for the class $[(A_k,\mu,\widetilde M_{kl})].$
This shows that the matrix algebra's structure on fibers
of considered bundles plays the same role that a
metric in usual theory of vector bundles
(recall that in usual theory of vector bundles
for finding of a stable inverse class for a given
class with a representative $\xi$
we must embed the vector bundle $\xi$ into a product
bundle and take its orthogonal complement there).

Let us consider the bundle space ${\cal P}^{k-1}{\mathop{\times}
\limits_{Gr_{k,kl}}}{\cal Q}^{l-1},$ $(k,l)=1$ (which was defined
on page \pageref{spaceBSP}). We also denote this space by
$\widetilde{Gr}_{k,kl}.$
It follows from Definition \ref{bdd} that there exists the
bundle map (over $Gr_{k,kl}$)
$$\lambda _{k,l} :\widetilde
{Gr}_{k,kl}\hookrightarrow Gr_{k,kl}\times \mathbb{C} P^{kl-1}$$
such that its restriction to any fiber (isomorphic to $\mathbb{C}P^{k-1}
\times \mathbb{C}P^{l-1}$) is the Segre's map $\mathbb{C}P^{k-1}\times
\mathbb{C}P^{l-1}\hookrightarrow \mathbb{C}P^{kl-1}.$ If we take
the composition of $\lambda_{k,l}$ with the projection onto
$\mathbb{C}P^{kl-1},$ then we obtain the map
$$\kappa_{k,l}:
\widetilde{Gr}_{k,kl}\hookrightarrow \mathbb{C}P^{kl-1}$$ such that
its restriction to any fiber (isomorphic to
$\mathbb{C}P^{k-1}\times \mathbb{C}P^{l-1}$)
is the classifying map for the exterior tensor product of the
canonical line bundles over
$\mathbb{C}P^{k-1}$ and $\mathbb{C}P^{l-1}.$

By $Gr$ denote the
$H$-group $\varinjlim_{(k,l)=1}Gr_{k,kl}$ (recall that
$\varinjlim_{(k,l)=1}Gr_{k,kl}\simeq BSU_{\otimes}$ as $H$-groups).
By $\widetilde{Gr}$ denote the direct limit
$\varinjlim_{(k,l)=1}\widetilde{Gr}_{k,kl}$
of the FBSP's bundle spaces over $Gr_{k,kl}$ and by
$\kappa$ denote the direct limit
$$\varinjlim_{(k,l)=1}\kappa_{k,l}\colon \widetilde{Gr}
\rightarrow \mathbb{C}P^\infty$$\label{kappa}
of the maps $\kappa_{k,l},\: (k,l)=1.$

It can be proved that the maps
$$\widetilde{\phi}_{kl,mn}\colon \widetilde{Gr}_{k,kl}\times \widetilde
{Gr}_{m,mn}\rightarrow \widetilde{Gr}_{km,klmn},\quad (km,ln)=1$$
(which correspond to the multiplication of FBSP,
consequently their restriction to any pair of fibers is the map
$$(\mathbb{C}P^{k-1}\times \mathbb{C}P^{l-1})\times(
\mathbb{C}P^{m-1}\times \mathbb{C}P^{n-1})\hookrightarrow
\mathbb{C}P^{km-1}\times \mathbb{C}P^{ln-1}{\rm )}$$
induce the structure of $H$-group on the $\widetilde{Gr}.$

Consider the set of pairs $(\xi_k,\eta_l)$
of $k$ and $l$-dimensional ($(k,l)=1$) complex vector bundles
(not necessarily with the structural group $\SU$)
over a finite $CW$-complex such that
$$\xi_k\otimes \eta_l\cong \zeta \otimes [kl],$$
where $\zeta$ is some geometric line bundle.
Define on this set the equivalence relation which is
analogous to the FAB (or FBSP) stable equivalence relation
(see Definition \ref{floating}). The tensor product of
pairs induces the group operation on the set of
equivalence classes.
It is easy to prove that the $H$-group $\widetilde{Gr}$
is the represented space for this
homotopic functor to category of Abelian groups.

\begin{proposition}
The map $\kappa \colon \widetilde{Gr}
\rightarrow \mathbb{C}P^\infty$ is
the homomorphism of $H$-groups. Moreover $\widetilde{Gr}$
is isomorphic to $\BSU_\otimes \times \mathbb{C}P^\infty
\times \mathbb{C}P^\infty$ as an $H$-group
(recall that $\mathbb{C}P^\infty$ is the $H$-group
with respect to the operation which corresponds to the tensor product of
complex line bundles; this $H$-structure
is induced by the Segre's maps $\mathbb{C}P^{k-1}
\times \mathbb{C}P^{l-1}\hookrightarrow \mathbb{C}P^{kl-1}$).
\end{proposition}

Consider also the ''one-half'' ${\cal P}^{k-1}$ of the bundle space
of the canonical FBSP ${\cal P}^{k-1}{\mathop{\times}\limits_{Gr_{k,kl}}}
{\cal Q}^{l-1}$ over $Gr_{k,kl}$ (recall that
${\cal P}^{k-1}{\mathop{\times}\limits_{Gr_{k,kl}}}{\cal Q}^{l-1}$
is the product over $Gr_{k,kl}$ of the bundles
${\cal P}^{k-1}$ and ${\cal Q}^{l-1}$
with fibers $\mathbb{C}P^{k-1}$ and $\mathbb{C}P^{l-1}$
respectively). By $\widehat{Gr}_{k,kl}$ \label{half}
we denote the bundle space ${\cal P}^{k-1}$
over $Gr_{k,kl}.$ Let $\widehat{Gr}$ be the
direct limit $\varinjlim_{(k,l)=1}\widehat{Gr}_{k,kl}.$
Note that the $H$-structure on $\widetilde{Gr}$
may be restricted to the subspace $\widehat{Gr}.$
It is not hard to prove that the space $\widehat{Gr}$
is the $H$-group with respect to the induced $H$-structure.

\begin{proposition}
\label{Hsp}
$H$-groups $\widehat{Gr}$ and $\BU_\otimes$
are isomorphic. Moreover they are isomorphic to the
direct product $\BSU_\otimes \times
\mathbb{C}P^\infty.$
\end{proposition}

\newpage

\section{Formal groups over Hopf algebras}

The aim of this section is to define some generalization
of the notion of formal group. More precisely, we consider the
analog of formal groups with coefficients belonging to
a Hopf algebra. We also study some example of a formal group
over a Hopf algebra, which generalizes the formal group
of geometric cobordisms.

Recently some important connections
between the Landweber-Novikov algebra and the formal group
of geometric cobordisms were established
\footnote{B. I. Botvinnik, V. M. Bukhshtaber, S. P. Novikov,
S. A. Yuzvinskii  Algebraicheskie aspekti teorii umnojenii
v kompleksnikh kobordizmakh // UMN. -- 2000. -- T.~55, ¹4.
S.~5--24.}.

Let $(H,\mu ,\eta ,\Delta ,\varepsilon, S)$ be a
(topological) Hopf algebra
over ring $R$ (where $\mu \colon H{\mathop{\widehat{\otimes}}\limits_R}H
\rightarrow H$ is multiplication, $\eta \colon R\rightarrow H$
is unit, $\Delta \colon H\rightarrow H
{\mathop{\widehat{\otimes}}\limits_R}H$
is diagonal (comultiplication), $\varepsilon \colon H\rightarrow
R$ is counit, and $S\colon H\rightarrow H$ is antipode).

\begin{definition}
\label{fgoha}
A formal series ${\frak F}(x\otimes 1,1\otimes x)\in
H{\mathop{\widehat{\otimes}}\limits_R}H[[x\otimes 1,1\otimes x]]$
is called {\it a formal group over the Hopf algebra}
$(H,\mu ,\eta ,\Delta ,\varepsilon ,S)$
if the following conditions hold:
\begin{itemize}
\item[1)](associativity)
$$((\id_H\otimes \Delta){\frak F})
(x\otimes 1\otimes 1,1\otimes {\frak F}(x\otimes 1,1\otimes x))
=$$
$$((\Delta \otimes \id_H){\frak F})({\frak F}(x\otimes 1,
1\otimes x)\otimes 1,1\otimes 1\otimes x);$$
\item[2)](unit) $$((\id_H\otimes \varepsilon){\frak F})(x\otimes 1,
0)=x\otimes 1,$$ $$((\varepsilon \otimes \id_H){\frak F})(0,1\otimes x)
=1\otimes x;$$
\item[3)](inverse element) there exists the series
$\Theta(x)\in H[[x]]$ such that
$$((\mu \circ (\id_H\otimes S)){\frak F})
(x,\Theta(x))=0=((\mu \circ (S\otimes \id_H))
{\frak F})(\Theta(x),x).$$
\end{itemize}
If for a formal group
${\frak F}(x\otimes 1,1\otimes x)$ over
a commutative and cocommutative Hopf algebra $H$
the equality ${\frak F}(x\otimes 1,1\otimes x)=
{\frak F}(1\otimes x,x\otimes 1)$ holds, then it is called
commutative.
Below we shall deal only with the commutative case.
\end{definition}

\begin{remark}
Note that a formal group ${\frak F}(x\otimes 1,1\otimes x)$
over Hopf algebra $H$ over ring $R$
defines the formal group (in the usual sense) $F(x\otimes
1,1\otimes x)\in R[[x\otimes 1,1\otimes x]]$ over the ring $R$
in the following way. By $F(x\otimes 1,1\otimes x)$
denote the series
$((\varepsilon \otimes
\varepsilon){\frak F})(x\otimes 1,1\otimes x).$
If we identify $R{\mathop{\otimes}\limits_R}R$ and $R,$
we may assume that $F(x\otimes 1,1\otimes x)\in
R[[x\otimes 1,1\otimes x]].$
Note that for any coalgebra $H$ the diagram
\begin{equation}
\label{coalg}
\begin{array}{ccc}
H & \stackrel{\Delta}{\rightarrow} & H{\mathop{\widehat{\otimes}}
\limits_R}H\\
\scriptstyle{\varepsilon}
\downarrow & & \quad \downarrow
\scriptstyle{\varepsilon \otimes \varepsilon}\\
R & \stackrel{\cong}{\rightarrow} & R\otimes R\\
\end{array}
\end{equation}
is commutative. Using (\ref{coalg}) and
condition 1) of Definition \ref{fgoha}, we get $F(x\otimes 1\otimes 1,
1\otimes F(x\otimes 1,1\otimes x))=F(F(x\otimes 1,1\otimes x)\otimes 1,
1\otimes 1\otimes x).$
Similarly the conditions $F(x\otimes 1,0)=x\otimes 1$ and
$F(0,1\otimes x)=1\otimes x$
may be verified.
It is well known,
that the existence of the inverse element (in the case of usual
formal groups) follows from
the proved conditions.
However this may be deduced from the condition 3) of Definition
\ref{fgoha} in the standard way.
Moreover, the inverse element
$\theta(x)$ in the formal group $F(x\otimes 1,1\otimes x)$
is equal to $(\varepsilon(\Theta))(x).$

Therefore we may consider the formal group
${\frak F}(x\otimes 1,1\otimes x)$ over Hopf algebra $H$
as an extension of the usual formal group $F(x\otimes 1,1\otimes x)$
by the Hopf algebra $H.$
\end{remark}

\begin{remark}
By definition, put
$\widetilde{\Delta}(x)={\frak F}(x\otimes 1,1\otimes x),\;
\widetilde{\varepsilon}(x)=0,\; \widetilde{S}(x)=\Theta(x)$
and $\widetilde{\Delta}\mid_H=\Delta, \;
\widetilde{\varepsilon}\mid_H=\varepsilon, \;
\widetilde{S}\mid_H=S.$
We claim that $(H[[x]],\widetilde{\mu},\widetilde{\eta},
\widetilde{\Delta},\widetilde{\varepsilon},\widetilde{S})$
is the Hopf algebra
(here $\widetilde{\mu},\; \widetilde{\eta}$ are evidently
extensions of $\mu,\; \eta$).
Indeed, the commutativity of the diagram
\begin{equation}
\begin{array}{ccc}
H[[x]] & \stackrel{\widetilde{\Delta}}{\rightarrow} &
H[[x]]{\mathop{\widehat{\otimes}}\limits_R}H[[x]] \\
{\scriptstyle \widetilde{\Delta}}\downarrow \; \: &&
\qquad \qquad \downarrow {\scriptstyle \id_{H[[x]]}
\otimes \widetilde{\Delta}} \\
H[[x]]{\mathop{\widehat{\otimes}}\limits_R}H[[x]] &
\stackrel{\widetilde{\Delta}\otimes \id_{H[[x]]}}
{\rightarrow} & H[[x]]{\mathop{\widehat{\otimes}}\limits_R}
H[[x]]{\mathop{\widehat{\otimes}}\limits_R}H[[x]] \\
\end{array}
\end{equation}
follows from the equations
$$(\id_{H[[x]]}\otimes \widetilde{\Delta})
(\frak F(x\otimes 1,1\otimes x))=
((\id_H \otimes \Delta){\frak F})(x\otimes 1\otimes 1,
1\otimes {\frak F}(x\otimes 1,1\otimes x))=$$
$$((\Delta \otimes \id_H){\frak F})({\frak F}(x\otimes 1,1\otimes x)
\otimes 1,1\otimes 1\otimes x)=(\widetilde{\Delta}\otimes \id_{H[[x]]})
({\frak F}(x\otimes 1,1\otimes x)).$$
The commutativity of the diagram
$$
\diagram
R{\mathop{\otimes}\limits_R}H[[x]] &
\qquad H[[x]]{\mathop{\widehat{\otimes}}\limits_R}H[[x]]
\qquad \lto_{\widetilde{\varepsilon}\otimes \id_{H[[x]]}}
\rto^{\id_{H[[x]]} \otimes \widetilde{\varepsilon}}&
H[[x]]{\mathop{\otimes}\limits_R}R \\
& H[[x]]\ulto^\cong \uto_{\widetilde{\Delta}} \urto_\cong & \\
\enddiagram
$$
follows from the equations
$$((\id_{H[[x]]}\otimes \widetilde{\varepsilon})\circ \widetilde{\Delta})
(x)=((\id_H \otimes \varepsilon){\frak F})(x\otimes 1,1\otimes
\widetilde{\varepsilon}(x))=x\otimes 1,$$
$$((\widetilde{\varepsilon}\otimes \id_{H[[x]]})\circ \widetilde{\Delta})
(x)=((\varepsilon \otimes \id_H){\frak F})
(\widetilde{\varepsilon}(x)\otimes 1,
1\otimes x)=1\otimes x.$$
The axiom of antipode
$$(\widetilde{\mu}\circ (\id_{H[[x]]}\otimes \widetilde{S})
\circ \widetilde{\Delta})(x)=
(\widetilde{\mu}\circ (\widetilde{S}\otimes \id_{H[[x]]})
\circ \widetilde{\Delta})(x)=
(\widetilde{\eta}\circ \widetilde{\varepsilon})(x)=0$$
follows from the condition 3) of Definition \ref{fgoha}.
\end{remark}

\begin{remark}
We may rewrite the conditions 1), 2), 3)
of Definition \ref{fgoha} in terms
of series ${\frak F}(x\otimes 1,1\otimes x)$ in the next way.
Let $$\sum_{i,j\geq 0}
A_{i,j}x^i\otimes x^j=$$
$$\sum_{i,j\geq 0}
(\sum_ka_{i,j}^k
\otimes b_{i,j}^k)x^i\otimes x^j\in
H{\mathop{\widehat{\otimes}}\limits_R}H [[x\otimes 1,1\otimes
x]]$$ be the series ${\frak F}(x\otimes 1,1\otimes x).$ Then
the condition 1) is equivalent to the following equality:
$$\sum_{i,j\geq 0}
(\sum_ka_{i,j}^k\otimes
\Delta(b_{i,j}^k))x^i\otimes {\frak F}(x\otimes 1,1\otimes
x)^j=$$
$$\sum_{i,j\geq 0}
(\sum_k\Delta(a_{i,j}^k)
\otimes b_{i,j}^k){\frak F}(x\otimes 1,1\otimes x)^i\otimes
x^j$$
The condition 2) is equivalent to
$$\sum_ka_{i,0}^k
\varepsilon (b_{i,0}^k)=0,\quad {\rm if}\quad i\neq 1,\quad
\sum_ka_{1,0}^k\varepsilon (b_{1,0}^k)=1,$$
$$\sum_k\varepsilon (a_{0,j}^k)b_{0,j}^k=0,\quad
{\rm if}\quad j\neq 1,\quad \sum_k\varepsilon (a_{0,1}^k)
b_{0,1}^k=1.$$
The condition 3) also may be rewritten in terms of series.
\end{remark}

Let us consider some examples of defined objects.

\begin{example}\label{ttt}(Trivial extension)
Let $F(x\otimes 1,1\otimes x)$
be a formal group (in the usual sense) over a ring $R,$
and $(H,\mu ,\eta ,\Delta ,\varepsilon ,S)$
be a Hopf algebra over the same ring $R.$
Then ${\frak F}(x\otimes 1,1\otimes x)=
((\eta \otimes \eta )F)(x\otimes 1,1\otimes x)\in
H{\mathop{\widehat{\otimes}}\limits_R}H[[x\otimes 1,1\otimes x]]$
is the formal group over the Hopf algebra $H$
(recall that we identify $R{\mathop{\otimes}\limits_R}R$
and $R$).
\end{example}

\begin{example} Now we construct a nontrivial
extension ${\frak F}(x\otimes 1,1\otimes x)$ of the formal group of
geometric cobordisms
$F(x\otimes 1,1\otimes x)\in \Omega_U^*(\pt)[[x\otimes 1,1\otimes x]]$
by the Hopf algebra
$\Omega_U^*(Gr).$ For this let us consider the map (see page
\pageref{half})
$\widehat{Gr}_{k,kl}\times \widehat{Gr}_{m,mn}
\stackrel{\widehat{\phi}_{kl,mn}}{\rightarrow}
\widehat{Gr}_{km,klmn},$ $(km,ln)=1.$
By $x\mid_{km,ln}$ denote the cobordism's class in
$\Omega_U^2(\widehat{Gr}_{km,klmn})$ such that
its restriction to every fiber of the bundle
\begin{equation}
\begin{array}{ccc}
\mathbb{C}P^{km-1} & \hookrightarrow & \widehat{Gr}_{km,klmn}\\
&& \downarrow \\
&& Gr_{km,klmn} \\
\end{array}
\end{equation}
is the standard generator in
$\Omega_U^2(\mathbb{C}P^{km-1}).$
Let $x\mid_{k,l}$ and $x\mid_{m,n}$ be analogously elements in
$\Omega_U^2(\widehat{Gr}_{k,l})$ and
$\Omega_U^2(\widehat{Gr}_{m,mn})$ respectively.
Then we obtain that $$\widehat{\phi}_{kl,mn}^*
(x\mid_{km,ln})=\sum_{0\leq i\leq k-1
\atop{0\leq j\leq m-1}}A_{i,j}
\mid_{kl,mn}(x\mid_{k,l})^i\otimes (x\mid_{m,n})^j,$$
where $A_{i,j}\mid_{kl,mn}\in \Omega_U^{2(1-i-j)}
(Gr_{k,kl}\times Gr_{m,mn}).$
Applying the functor of unitary cobordisms to the
following injective system
of the spaces and their maps
\begin{equation}
\begin{array}{ccc}
\widehat{Gr}_{p,pq}\times
\widehat{Gr}_{t,tu} &
\stackrel{\widehat{\phi}_{pq,tu}}{\rightarrow} &
\widehat{Gr}_{pt,pqtu}\\
\uparrow \quad && \uparrow \qquad \\
\widehat{Gr}_{k,kl}\times
\widehat{Gr}_{m,mn} &
\stackrel{\widehat{\phi}_{kl,mn}}{\rightarrow} &
\widehat{Gr}_{km,klmn}.\\
\end{array}
\end{equation}
(under the conditions $k\mid p,\quad l\mid q,\quad
m\mid t,\quad n\mid u,\quad$and
$(pt,qu)=1$),
we obtain the formal series $${\frak F}(x\otimes 1,1\otimes x)=$$
$$\sum_{i,j\geq 0}A_{i,j}x^i\otimes x^j\quad \in \Omega_U^*
(Gr){\mathop{\widehat{\otimes}}
\limits_{\Omega_U^*(\pt)}}\Omega_U^*(Gr)[[x\otimes 1,
1\otimes x]]$$
such that $i^*_{kl}A_{i,j}=A_{i,j}\mid_{k,l}$ for injection
$i_{kl}\colon Gr_{k,kl}\hookrightarrow Gr$ (for every pair $\{
k,l\}$ such that $(k,l)=1$).

By $R$ and $H$ denote the ring
$\Omega_U^*(\pt)$ and the Hopf algebra $\Omega_U^*(Gr)$ (over
the ring $\Omega_U^*(\pt)$) respectively
(recall that we consider the space $Gr$ with
the $H$-group structure, induced by the multiplication of
FBSP).

\begin{proposition}
\label{gfggc}
The series
${\frak F}(x\otimes 1,1\otimes x)$
is the formal group over the Hopf algebra $H.$
\end{proposition}
{\raggedright {\it Proof}.}
To prove $((\id_H\otimes \Delta){\frak F})(x\otimes 1\otimes 1,1\otimes
{\frak F}(x\otimes 1,1\otimes x))=((\Delta \otimes \id_H){\frak F})
({\frak F}(x\otimes 1,1\otimes x)\otimes 1,1\otimes 1\otimes x),$
we need the following commutative diagram ($(kmt,lnu)=1$):
\begin{equation}
\begin{array}{ccc}
\widehat{Gr}_{k,kl}\times \widehat{Gr}_{m,mn}
\times \widehat{Gr}_{t,tu} & \rightarrow &
\widehat{Gr}_{k,kl}\times
\widehat{Gr}_{mt,mntu}\\
\downarrow && \downarrow \\
\widehat{Gr}_{km,klmn}\times \widehat{Gr}_{t,tu} &
\rightarrow & \widehat{Gr}_{kmt,klmntu}. \\
\end{array}
\end{equation}

To prove $((\id_H\otimes \varepsilon){\frak F})(x\otimes 1,0)=
x\otimes 1,$
we need the following commutative diagram ($(km,ln)=1$):
\begin{equation}
\begin{array}{ccc}
\widehat{Gr}_{k,kl}\times \widehat{Gr}_{m,mn} & \rightarrow &
\widehat{Gr}_{km,klmn} \\
\uparrow && \uparrow \\
\widehat{Gr}_{k,kl}\times \mathbb{C}P^{m-1} &
\leftarrow & \widehat{Gr}_{k,kl}\times \{ \pt \}, \\
\end{array}
\end{equation}
where right-hand vertical arrow is the standard inclusion.

To prove $((\mu \circ (\id_H\otimes S)){\frak F})(x,\Theta
(x))=0,$ let us construct the fiber map $\widehat{\nu}\colon
\widehat{Gr}\rightarrow \widehat{Gr}$ such that the
following two conditions are satisfied:
\begin{itemize}
\item[1)] the restriction of $\widehat{\nu}$
to any fiber ($\cong \mathbb{C}P^\infty$)
is the inversion in the $H$-group $\mathbb{C}P^\infty$;
\item[2)] $\widehat{\nu}$ covers the $\nu \colon Gr\rightarrow Gr$
(where $\nu$ is the inversion in the $H$-group $Gr$).
\end{itemize}

Let us remember that
${\cal P}^{k-1}{\mathop{\times}\limits_{Gr_{k,kl}}}{\cal Q}^{l-1}$
is the canonical FBSP over $Gr_{k,kl}$ and we have denoted
by $\widehat{Gr}_{k,kl}$ the bundle space ${\cal P}^{k-1}.$
Let $\widehat{Gr}'_{k,kl}$ ($\widehat{Gr}'$) be the bundle space of the
''second half'' ${\cal Q}^{l-1}$ of the canonical FBSP over
$Gr_{k,kl}$ ($\varinjlim_{(k,l)=1}\widehat{Gr}'_{k,kl}$
respectively).

First note that there exists the fiber isomorphism
$\widehat{\nu}'_{k,l}\colon \widehat{Gr}_{k,kl}
\rightarrow \widehat{Gr}'_{l,lk}$
that covers the inverse map $\nu_{k,l}\colon Gr_{k,kl}\rightarrow
Gr_{l,lk}$ (in other words, the map $\nu_{k,l}$ takes each
subalgebra $A_k\cong M_k(\mathbb{C})$ in the $M_{kl}(\mathbb{C})$
to its centralizer $\Z_{M_{kl}(\mathbb{C})}(A_k)\cong M_l(\mathbb{C})$
in the $M_{kl}(\mathbb{C})$). Let $c_{l,k}\colon
\widehat{Gr}'_{l,lk}\rightarrow \widehat{Gr}'_{l,lk}$
be the fiber map such that the following two conditions
are satisfied:
\begin{itemize}
\item[1)] $c_{l,k}$ covers the identity mapping of the base $Gr_{l,lk};$
\item[2)] the restriction of $c_{l,k}$ to any fiber
$\cong \mathbb{C}P^{k-1}$ is the complex conjugation.
\end{itemize}
Let $\widehat{\nu}_{k,l}\colon \widehat{Gr}_{k,kl}\rightarrow
\widehat{Gr}'_{l,lk}$ be the composition
$c_{l,k}\circ \widehat{\nu}'_{k,l}.$
It is easy to prove that the map
$\widehat{\nu}=\varinjlim_{(k,l)=1}\widehat{\nu}_{k,l}
\colon \varinjlim_{(k,l)=1}\widehat{Gr}_{k,kl}
\rightarrow \varinjlim_{(l,k)=1}\widehat{Gr}'_{l,lk}$
is required. In particular, there exists the fiber isomorphism
between $\widehat{Gr}$ and $\widehat{Gr}'.$

The map $\widehat{\nu}$ defines (by the same way, as
$\widehat{\phi}$ in the beginning of the example) the formal series
$\Theta(x)\in H[[x]]$ (note that
$\varepsilon(\Theta)(x)=
\theta(x),$ where $\theta(x)\in R[[x]]$ is the inverse element in
the group of geometric cobordisms).

Now we claim that
$((\mu \circ (\id_H\otimes S)){\frak F})(x,\Theta(x))=0.$
Indeed, this follows from the next commutative diagram:
\begin{equation}
\begin{array}{ccccccc}
\widehat{Gr} & \stackrel{\diag}{\rightarrow} &
\widehat{Gr}\times \widehat{Gr} &
\stackrel{\id \times \widehat{\nu}}{\rightarrow} &
\widehat{Gr}\times \widehat{Gr} & \stackrel{\widehat{\phi}}
{\rightarrow} & \widehat{Gr} \\
\downarrow && \downarrow && \downarrow && \downarrow \\
Gr & \stackrel{\diag}{\rightarrow} & Gr\times Gr &
\stackrel{\id \times \nu}{\rightarrow} & Gr\times Gr &
\stackrel{\phi}{\rightarrow} & Gr \\
\end{array}
\end{equation}
(we see that the composition
$\widehat{\phi} \circ (\id \times \widehat{\nu})\circ \diag$
is homotopic (in class of fiber homotopies) to the
map $\widehat{Gr}\rightarrow \pt \in \widehat{Gr}$).
$\square$

\smallskip

In section 1 we considered the map $\kappa \colon
\widetilde{Gr}\rightarrow \mathbb{C}P^\infty.$
It is the direct limit of the fiber maps (see p. \pageref{kappa})
\begin{equation}
\begin{array}{ccc}
\widetilde{Gr}_{k,kl} & \stackrel{\kappa_{k,l}}{\rightarrow} &
\mathbb{C}P^{kl-1} \\
\downarrow && \downarrow \\
Gr_{k,kl} & \rightarrow & \pt. \\
\end{array}
\end{equation}
It defines (in the same way, as $\widehat{\phi}$ and
$\widehat{\nu}$ above) the formal
series
$${\frak G}(x,y)=
\sum_{i,j\geq 0}B_{i,j}x^iy^j\; \in H[[x,y]].$$

\begin{proposition}
${\frak G}(x,y)=
((\mu \circ (\id_H\otimes S)){\frak F})(x,y),$
i. e. $$B_{i,j}=\sum_ka_{i,j}^kS(b_{i,j}^k).$$
\end{proposition}
{\raggedright {\it Proof}.}
Recall that in the proof of Proposition \ref{gfggc}
the fiber maps $\widehat{\nu}'_{k,l}\colon \widehat{Gr}_{k,kl}
\rightarrow \widehat{Gr}'_{l,lk}$ were defined.
By $\widehat{\nu}'$ denote the direct limit $\varinjlim_{(k,l)=1}
\widehat{\nu}'_{k,l}\colon \widehat{Gr}\rightarrow \widehat{Gr}.$
Note that $\widehat{\nu}'$ covers the inversion
$\nu \colon Gr\rightarrow Gr$ in the $H$-group $Gr.$

Now the proof follows from the next composition of the bundle maps:
\begin{equation}
\begin{array}{ccccccc}
\widetilde{Gr} & \rightarrow &
\widehat{Gr}\times \widehat{Gr} & \stackrel{\id \times \widehat{\nu}'}
{\rightarrow} &
\widehat{Gr}\times
\widehat{Gr} & \stackrel{\widehat{\phi}}{\rightarrow} & \widehat{Gr} \\
\downarrow && \downarrow && \downarrow && \downarrow \\
Gr & \stackrel{diag}{\rightarrow} & Gr\times Gr & \stackrel{\id \times
\nu}{\rightarrow} & Gr\times Gr & \stackrel{\phi}{\rightarrow} & Gr. \\
\end{array}
\end{equation}
We see that the upper composition in fact is the map $\widetilde{Gr}
\rightarrow \mathbb{C}P^\infty$ and it coincides with the map $\kappa.$
Let $y$ be $\widehat{\nu}'{^*}(x).$ The upper composition gives
$x\mapsto {\frak F}(x\otimes 1,1\otimes x)\mapsto
((\id_H\otimes S){\frak F})(x\otimes 1,1\otimes y)\mapsto
((\mu \circ (\id_H\otimes S)){\frak F})(x\otimes 1,1\otimes y).$
Without loss of sense we may write $x$ and $y$ instead of
$x\otimes 1$ and $1\otimes y$ respectively. $\square$

The series ${\frak G}(x,y)$ has the following
interesting property.
\begin{proposition}
$$(\Delta{\frak G})({\frak F}(x\otimes 1,1\otimes x),
((S\otimes S){\frak F})(y\otimes 1,1\otimes y))=$$
$$F({\frak G}(x,y)\otimes 1,1\otimes{\frak G}(x,y)),$$
where $F(x,y)\; \in R[[x,y]]$ is the formal group
of geometric cobordisms.
\end{proposition}
{\raggedright {\it Proof}.}
We give two variants of the proof.

1).''Topological proof'' follows from the commutative diagram
\begin{equation}
\begin{array}{ccc}
\mathbb{C}P^{kl-1}\times \mathbb{C}P^{mn-1} & \rightarrow &
\mathbb{C}P^{klmn-1} \\
\uparrow \; \: && \uparrow \\
\widetilde{Gr}_{k,kl}\times \widetilde{Gr}_{m,mn} & \rightarrow &
\widetilde{Gr}_{km,klmn} \\
\end{array}
\end{equation}
($(km,ln)=1$) combining with the decomposition of the map $\kappa,$
which was obtained in previous proof.

2). By $\widetilde{S}'$ denote the homomorphism
$\widehat{\nu}'{^*}
\colon H[[x]]\rightarrow H[[y]]$ (recall that $\widehat{\nu}'{^*}\mid_H=S
\colon H\rightarrow H,$ where $S$ is the antipode).
Let us consider the following composition of homomorphisms
of Hopf algebras:
$$H[[x]]\stackrel{\widetilde{\Delta}}{\rightarrow}
H[[x]]{\mathop{\widehat{\otimes}}\limits_R}H[[x]]
\stackrel{\id \otimes \widetilde{S}'}{\rightarrow}
H[[x]]{\mathop{\widehat{\otimes}}\limits_R}H[[y]]
\stackrel{(\mu)}{\rightarrow}H[[x,y]],$$
where $(\mu)$ is the homomorphism, induced by multiplication
$$\mu \colon H{\mathop{\widehat{\otimes}}\limits_R}H\rightarrow H.$$
It follows from the axiom of antipode $$\mu \circ (\id_H\otimes S)\circ
\Delta =\eta \circ \varepsilon$$ that $$(\mu)\circ (\id_{H[[x]]}\otimes
\widetilde{S}')\circ \widetilde{\Delta}\mid_H=\eta \circ \varepsilon.$$
Hence there exists the homomorphism of Hopf algebras
$$(\eta)\colon R[[x]]\rightarrow H[[x,y]]$$
such that the following diagram
$$
\diagram
H[[x]] \rto^{{\widetilde{\Delta}}\quad}
\drto_{(\varepsilon)} & H[[x]]{\mathop{\widehat{\otimes}}\limits_R}H[[x]]
\rto^{\id \otimes \widetilde{S}'} &
H[[x]]{\mathop{\widehat{\otimes}}\limits_R}H[[y]] \dto^{(\mu)} \\
& R[[x]] \rto^{(\eta)\qquad} & H[[x,y]] \\
\enddiagram
$$
is commutative (here $(\varepsilon)$ is the homomorphism,
induced by $\varepsilon$). Note that $$(\eta)(x)=
{\frak G}(x,y).$$
This completes the proof that $$\Delta_{R[[x]]}(x)=F(x\otimes 1,
1\otimes x),$$ where $F(x\otimes 1,1\otimes x)\in R[[x\otimes
1,1\otimes x]]$ is the formal group of geometric cobordisms. $\square$

It is very important that we consider the maps
$\widehat{\phi},\; \widehat{\nu},$ and $\kappa$ as {\it fiber} maps
in this example. Otherwise instead of
${\frak F}(x\otimes 1,1\otimes x)$ we obtain the usual formal group
of geometric cobordisms because the $H$-space $\widehat{Gr}$
is isomorphic to the $H$-space $BSU_\otimes \times \mathbb{C}P^\infty$
(see Proposition \ref{Hsp}).
\end{example}

It is well known\footnote{Quillen D. On the formal group low
of unoriented and complex cobordism theory, Bull. Amer. Math. Soc.,
75:6 (1969), 1293--1298.}, that the formal group
of geometric cobordisms is the universal
formal group.
\begin{conjecture}
\nonumber
The formal group ${\frak F}(x\otimes 1,1\otimes x)$
is the universal object in the category of formal groups
over a (topological) Hopf algebras.
\end{conjecture}

Let $R'$ be a ring and $F'(x\otimes 1,1\otimes x)$
be a formal group over $R'.$ Note that we may
consider the $R'$ as the Hopf algebra over $R'$
with respect to the $\Delta_{R'}\colon R'\cong R'{\mathop{\otimes}
\limits_{R'}}R',$ $\eta_{R'}=\varepsilon_{R'}=S_{R'}=\id_{R'}
\colon R'\rightarrow R'.$
If $\chi \colon H\rightarrow R'$ is a homomorphism of the
Hopf algebras from $(H,\mu,\eta,\Delta,\varepsilon,S)$ to
$(R',\mu_{R'},\eta_{R'},\Delta_{R'},\varepsilon_{R'},S_{R'}),$
then $\chi = (\chi \circ \eta)\circ \varepsilon =\chi \mid _R
\circ \varepsilon.$ Hence there exists the natural bijection
$\Hom_{Hopf\: alg.}(H,R')\leftrightarrow \Hom_{Ring}(R,R').$
Therefore the Conjecture
implies the universal property of the
formal group of geometric cobordisms.

\smallskip

The author is grateful to professor E. V. Troitsky
for constant attention to this work, and to professors
V. M. Manuilov and A. S. Mishchenko for useful discussions.

\end{document}